\documentclass[11pt]{article} 
\usepackage{amsmath,amssymb}

\usepackage{texdraw}
\everytexdraw{\drawdim pt \linewd 0.5 \textref h:C v:C \setunitscale 1
  \arrowheadtype t:F \arrowheadsize l:5 w:3
}

\newcommand{\onedot}{\bsegment \move (0 0) \fcir f:0 r:2 \esegment}
\newcommand{\smalldot}{\bsegment \fcir f:0 r:1.5 \esegment}

\usepackage[all]{xy}
\SelectTips{cm}{}
\newdir{ >}{{}*!/-7pt/@{>}}
\entrymodifiers={+!!<0pt,\fontdimen22\textfont2>}

\newcommand{\drpullback}[1][dr]{\save*!/#1-1.2pc/#1:(-1,1)@^{|-}\restore}

\newcommand{\G}{\mathbb{G}}
\newcommand{\N}{\mathbb{N}}

\newcommand{\Hom}{\operatorname{Hom}}
\newcommand{\Aut}{\operatorname{Aut}}
\DeclareMathOperator*{\colim}{colim}

\newcommand{\op}{^{\text{{\rm{op}}}}}
\providecommand{\kat}[1]{\text{\textbf{\textsl{#1}}}}
\newcommand{\Set}{\kat{Set}}
\newcommand{\Cat}{\kat{Cat}}
\newcommand{\sSet}{\kat{sSet}}
\newcommand{\Sh}{\kat{Sh}}
\newcommand{\PrSh}{\kat{PrSh}}
\newcommand{\Gr}{\kat{Gr}}
\newcommand{\elGr}{\kat{elGr}}
\newcommand{\triv}{\star}

\newcommand{\comma}{\raisebox{1pt}{$\downarrow$}}
\newcommand{\into}{\hookrightarrow}
\newcommand{\isopil}{\stackrel{\raisebox{0.1ex}[0ex][0ex]{\(\sim\)}}%
			{\raisebox{-0.15ex}[0.28ex]{\(\rightarrow\)}}}

\newcommand{\ov}{\overline}

\makeatletter

\def\abstractsize{\fontsize{\@viiipt}{8pt}\selectfont}

\normalsize                                
\newdimen\@bls                             
\@bls=\baselineskip                        
\topmargin \z@           
\headheight  \z@         
\headsep     \z@         
\bigskipamount=\@bls \@plus 0.3\@bls \@minus 0.3\@bls 
\advance\textheight\topskip  
\columnseprule \z@           
\skip\footins 12\p@ \@plus  8\p@      
\@fptop \z@ \@plus 1fil  
\@fpbot \z@ \@plus 1fil  
\@dblfptop \z@ \@plus 1fil 
\@dblfpbot \z@ \@plus 1fil 
\newskip\eqntopsep                   
\newdimen\eqnarraycolsep             
\def\section{\@startsection{section}{1}{\z@}%
        {-3ex plus -1ex minus -.2ex}%
        {2ex plus .2ex}{\large\bfseries}}

\def\lb@section{\thesection.\half@em}
\def\lb@empty@section{\thesection}

\def\head@style{\interlinepenalty\@M \hyphenpenalty\@M
  \exhyphenpenalty\@M \rightskip 0pt plus 0.5\hsize \relax}

\tabbingsep \labelsep  
\skip\@mpfootins = 6\p@ \@plus 2\p@  
\def\@captionsize{\small}

\long\def\@maketablecaption#1#2{\@captionsize
  \hbox to \hsize{\parbox[t]{\hsize}{\begin{center}#1 \\ #2\end{center}}}}
\long\def\@makefigurecaption#1#2{\@captionsize
  \vskip 8\p@
  \setbox\@tempboxa\hbox{#1. #2}
  \ifdim \wd\@tempboxa >\hsize              
    \unhbox\@tempboxa\par                   
  \else                                     
    \hbox to\hsize{\hfil\box\@tempboxa\hfil}
  \fi}
\def\conttablecaption{\par \begingroup \@parboxrestore \normalsize
  \@makecaption{\fnum@table\,---\,continued}{}\par
  \vskip-1pc \endgroup}
\def\contfigurecaption{\vskip-1pc \par \begingroup \@parboxrestore \normalsize
  \@makecaption{\fnum@figure\,---\,continued}{}\par
  \endgroup}

\def\thetable{\@arabic\c@table}
\def\fps@table{tbp}
\def\ftype@table{2}
\def\ext@table{lot}
\def\fnum@table{\tablename~\thetable}
\def\table{%
\let\@makecaption\@maketablecaption
\small
  \let\footnoterule\relax
  \let\contcaption\conttablecaption \@float{table}}
\def\endtable{\addvspace{.2pc plus .2pc minus .1pc}\end@float}
\@namedef{table*}{%
\let\@makecaption\@maketablecaption
\small
  \let\footnoterule\relax
  \let\contcaption\conttablecaption \@dblfloat{table}}
\@namedef{endtable*}{\end@dblfloat}

\makeatother

\begin{document}

\begin{center}
  
  QPL 2009

  \vspace{120pt}
  
  \textbf{\LARGE Feynman graphs, and nerve theorem for \\[1pt]
  compact symmetric multicategories \\[3pt]
  (extended abstract)}
  
  \vspace{16pt}
  
  {\Large Andr\'e Joyal}
  \footnote{Supported by the NSERC of Canada}
  
  \vspace{4pt}
  
  {\abstractsize {\em 
  {D\'epartement des math\'ematiques\\ Universit\'e du Qu\'ebec \`a 
    Montr\'eal\\[-5pt] Canada}
    }}
    
    \vspace{12pt}
    
  {\Large Joachim Kock}
  \footnote{Email: \texttt{kock@mat.uab.cat}}
  \footnote{Supported by research grants MTM2006-11391 and
  MTM2007-63277 of Spain.}
  
  \vspace{4pt}
  
     {\abstractsize {\em 
Departament de matem\`atiques\\ Universitat Aut\`onoma de Barcelona\\[-5pt]
    Spain} 
 }
  
\end{center}
 
 \vspace{12pt}
 
 \hrule
 
 \vspace{8pt}
 
 \abstractsize
 
 \noindent
 {\bf Abstract}

 \vspace{5pt}
 
 \noindent
  We describe a category of Feynman graphs and show how it relates to compact
  symmetric multicategories (coloured modular operads) just as linear orders
  relate to categories and rooted trees relate to multicategories.  More
  specifically we obtain the following nerve theorem: compact symmetric
  multicategories can be characterised as presheaves on the category of Feynman
  graphs subject to a Segal condition.  
  This text is a write-up of the second-named author's QPL6 talk; 
  a more detailed account of this material will
  appear elsewhere~\cite{Joyal-Kock:FeynmanInPrep}.

  \vspace{8pt}
  
  \noindent
{\em Keywords:}
  Feynman graph, multicategory, modular operad, nerve theorem, monad.

  \vspace{10pt}

  \hrule

  \normalsize
  
\vspace{8pt}

\section{Introduction}
\label{intro}

The graphical calculus of string diagrams is an important ingredient in
many abstract approaches to quantum mechanics and quantum informatics, as
well exemplified in this volume.
The objects of a category are pictured as strings,
and the arrows as dots.  An arrow is thought of as an operation, and in a plain
category each operation has precisely one input (source) and one output
(target).  Arrows can be composed if arranged in a sequence such that
the output of one operation matches the input of the next.
Table~\ref{table} below should illustrate the passage from
categories to compact symmetric multicategories.

\begin{table}[h!]
  \begin{center}
  \newcommand{\inlineDotlessTree}{%
\raisebox{4pt}{
\begin{texdraw} \linewd 0.5 \bsegment
    \move (0 0) \lvec (0 12) \move (1 0)
  \esegment \end{texdraw} } }

    \begin{tabular}{ l|| c|c| c l }
    \hline
    structure & objects & operations & pastings &
    \\ 
    \hline
    \raisebox{6pt}{categories \ } & \inlineDotlessTree
&
\begin{texdraw} \rlvec (0 6) \smalldot \rlvec (0 6) \move ( 0 -3)\end{texdraw}
& 
\begin{texdraw} \rlvec (0 5) \smalldot \rlvec (0 6) \smalldot \rlvec (0 6) 
  \smalldot \rlvec (0 4) \move (0 26) \end{texdraw}  & \raisebox{6pt}{linear orders}
  \\ 
    \raisebox{6pt}{multicategories \ } & \inlineDotlessTree& \begin{texdraw}
    \move (-5 0) \rlvec (0 5) \smalldot
    \move (10 0) \rlvec (0 5) \smalldot \rlvec (0 6)
    \move (25 0) \rlvec (0 5) \smalldot \rlvec (-3 6) \rmove (3 -6) \rlvec (3 6)
    \move (40 0) \rlvec (0 5) \smalldot \rlvec (-4 6) \rmove (4 -6) 
    \rlvec (0 7) \rmove (0 -7) \rlvec (4 6) \move ( 0 -3)
\end{texdraw}
&
\begin{texdraw}
  \lvec (0 6) \smalldot \lvec (-4 12) \smalldot \lvec (-6 20) \move (-4 12) 
  \lvec (-2 20) \move (0 6) \lvec (2 19) \move (0 6) \lvec (5 13) \smalldot 
  \move (0 26)
  \end{texdraw} & \parbox[b]{12ex}{planar \\[-2pt] rooted trees}
 \\
 \parbox[b]{20ex}{symmetric \\[-2pt] multicategories}
     & \inlineDotlessTree & \begin{texdraw}
    \move (-5 0) \rlvec (0 5) \smalldot
    \move (10 0) \rlvec (0 5) \smalldot \rlvec (0 6)
    \move (25 0) \rlvec (0 5) \smalldot \rlvec (-3 6) \rmove (3 -6) \rlvec (3 6)
    \move (45 0) \rlvec (0 5) \smalldot \rlvec (-4 6) \rmove (4 -6) 
    \rlvec (0 7) \rmove (0 -7) \rlvec (4 6)
    \htext (31 6) {\tiny $\mathfrak S\!_2$}
    \htext (52 5) {\tiny $\mathfrak S\!_3$}
    \move ( 0 -3)
\end{texdraw} 
&
\begin{texdraw}
  \lvec (0 6) \smalldot \lvec (-4 12) \smalldot \lvec (-6 20) \move (-4 12) 
  \lvec (-2 20) \move (0 6) \lvec (2 19) \move (0 6) \lvec (5 13) \smalldot
    \move (0 26)
  \end{texdraw}
  & \raisebox{6pt}{rooted trees}
\\ 
 \parbox[b]{20ex}{cyclic symmetric \\[-2pt] multicategories}
& \inlineDotlessTree & 
    \begin{texdraw}
\move (20 6) \smalldot
\move (34 6) \smalldot \lvec (34 13)
\move (48 6) \smalldot \lvec (53 12) \move (48 6) \lvec (53 0)
\move (70 6) \smalldot \lvec (65 12) \move (70 6) \lvec (65 0) \move (70 6) 
\lvec (77 6)
\move (90 6) \smalldot \move (84 11) \lvec (96 1)
\move (84 1) \lvec (96 11)
    \htext (56 6) {\tiny $\mathfrak S\!_2$}
    \htext (75 10) {\tiny $\mathfrak S\!_3$}
    \htext (100 6) {\tiny $\mathfrak S\!_4$}
\move ( 20 -3)
\end{texdraw}  & 
\begin{texdraw} \lvec (7 2) \smalldot \rlvec (6 5) \smalldot \rlvec (-3 6)\smalldot
  \move (13 7) \rlvec (7 3) \smalldot \rlvec (7 3) \rmove ( -7 -3) \rlvec (5 -4) 
  \move (13 7) \rlvec (5 -5) \smalldot \rlvec (0 -7) \rmove (0 7) \rlvec (6 
  -2)
  \move (7 2) \rlvec (-3 7) \move (7 2) \rlvec (2 -7) 
  \move (0 20)
   \end{texdraw}
   & \raisebox{6pt}{trees} \\ 
    \parbox[b]{20ex}{compact symmetric \\[-2pt] multicategories}
 &\inlineDotlessTree& 
    \begin{texdraw}
\move (20 6) \smalldot
\move (34 6) \smalldot \lvec (34 13)
\move (48 6) \smalldot \lvec (53 12) \move (48 6) \lvec (53 0)
\move (70 6) \smalldot \lvec (65 12) \move (70 6) \lvec (65 0) \move (70 6) 
\lvec (77 6)
\move (90 6) \smalldot \move (84 11) \lvec (96 1)
\move (84 1) \lvec (96 11)
    \htext (56 6) {\tiny $\mathfrak S\!_2$}
    \htext (75 10) {\tiny $\mathfrak S\!_3$}
    \htext (100 6) {\tiny $\mathfrak S\!_4$}
    \move ( 20 -3)
\end{texdraw} 
&
\begin{texdraw} \clvec (-12 -8)(-12 8)(0 0) \smalldot
  \clvec (0 8)(2 10)(10 10) \smalldot \move (0 0)
  \clvec (8 0)(10 2)(10 10)  \move (0 0) \lvec (0 -8)
  \move (0 0) \lvec (10 10) 
  \move (10 10) \lvec (17 3) \smalldot \lvec (11 -4) \smalldot
  \lvec (0 0) \lvec (-6 9)
  \move (17 3) \lvec (22 -4) \smalldot
    \move (0 18)
  \end{texdraw}
  & \raisebox{6pt}{Feynman graphs}
  \\
    \hline
    \end{tabular}
  \end{center}
  \caption{}
    \label{table}
\end{table}
%

{\em Multicategories} generalise categories by allowing a finite list of inputs
while still insisting on exactly one output for each operation.  The operations
are represented as planar rooted trees with precisely one inner vertex, and
composition of operations produces a single such one-vertex tree from any formal
configuration of matching operations, i.e.~turns an arbitrary tree into a
one-vertex tree.  (Monoidal categories are
a special kind of multicategories: they are the representable ones 
\cite{Leinster:0305049}.)
Abandonning the linear order on the set of inputs, we arrive
at the notion of {\em symmetric multicategory}, also known as {\em coloured
operad}.
Multicategories were introduced in 1969 by
Lambek~\cite{Lambek:deductiveII} to model sequent calculus, while operads,
the one-object symmetric analogue, often enriched over vector spaces or
topological spaces, were discovered at about the same time in loop
space theory \cite{Boardman-Vogt:LNM347}, \cite{May:LNM271}.

Giving up the distinction between input and output, calling the loose edges {\em
ports}, we arrive at a many-object version of what in topology is called {\em
cyclic operad}~\cite{Getzler-Kapranov:cyclic}.  Without a notion of input and
output, we impose instead an involution on the set of objects: each object has a
dual object.  The operations are now non-rooted (non-planar) one-vertex trees,
each edge of which is decorated by some object.  These operations can be
connected to each other by clutching a port of one operation to a matching
port of another operation (i.e.~the two decorating objects are required to be
dual), and the configurations of formal composites are non-rooted (non-planar)
trees.

Finally, by allowing also clutching of two ports of one and the same operation, we arrive
at the object of study of this work: the {\em compact symmetric multicategories}.
This is essentially a many-object version of the {\em modular operads} of
Getzler and Kapranov~\cite{Getzler-Kapranov:9408}, introduced to describe the
algebraic structure of the moduli space of curves in algebraic geometry, and in
topological and conformal field theory.  The configurations of formal
composites are now general connected graphs,  more precisely
what we call {\em Feynman graphs}: they are (non-directed) graphs,
allowed to have multiple edges and loops, as well as open edges.  We shall only use
connected Feynman graphs.

\pagebreak

Hence compact symmetric multicategories relate to Feynman graphs as categories
relate to linear orders.  It is the goal of this work to make this statement
precise.



\section{Nerve theorem for categories}
\label{Sec:classical}

We start by reviewing the classical nerve theorem,
following 
Berger~\cite{Berger:Adv}, Leinster~\cite{Leinster:CT04}, and
Weber~\cite{Weber:TAC18}.
Let $\Delta$ denote the category of finite
(non-empty) linearly ordered sets and monotone maps, and
recall that a simplicial set is a presheaf on $\Delta$, i.e.~a
functor $\Delta\op\to \Set$.  The {\em nerve} of a small
category $C$ is the simplicial set $NC$ whose $k$-simplicies are the chains of
 $k$ composable arrows in $C$.  
More conceptually, via the natural embedding $i:\Delta \into \Cat$
(interpreting an ordered set as a small category), $NC$ is simply the
presheaf
\begin{align}
  \Delta\op \  & \longrightarrow \ \Set \notag \\
  [k] \ & \longmapsto \ \Hom_{\Cat}(i([k]),C). \label{classicalnerve}
\end{align}
The {\em nerve theorem}, first observed by Grothendieck, asserts that
$N: \Cat \to \sSet$
is a fully faithful functor and that its essential image consists of those simplicial 
sets $X$
that satisfy the {\em Segal condition}: for each $k \geq 1$ the natural map
$$
X_k \longrightarrow X_1 \times_{X_0} \dots \times_{X_0} X_1
$$
is a bijection.  (The fibre product expresses composability: the target of one
arrow equals the source of the next.)  The nerve functor plays a fundamental
role to link category theory to topology, remembering that simplicial sets is
the most important combinatorial model for homotopy theory.

The main result of the present work is a nerve theorem for compact symmetric
multicategories which is a direct generalisation of the classical nerve theorem:
we characterise compact symmetric multicategories as presheaves on a category of
Feynman graphs satisfying a certain Segal condition. 
There are completely analogous nerve theorems for the other rows of
Table~\ref{table}.  (A nerve theorem in a slightly different spirit was obtained
for symmetric multicategories by Moerdijk and
Weiss~\cite{Moerdijk-Weiss:0701293}.)
The proof has two 
ingredients: one is to identify the correct category of graphs; the second is
an application of the abstract machinery developed by 
Berger~\cite{Berger:Adv}, Leinster~\cite{Leinster:CT04}, and
Weber~\cite{Weber:TAC18}.

A short review of the
classical case will be helpful.
A small category has
an underlying directed graph (see \cite{MacLane:categories}, 
Ch.II). A directed graph can be seen as a
presheaf on the category $\G =\{ 0\rightrightarrows 1 \}$.
The forgetful
functor from $\Cat$ to $\PrSh(\G)$ has a left adjoint, the {\em free-category}
functor (\cite{MacLane:categories}, Ch.II.7): 
the free category on a directed graph $G$ has the vertices of $G$ as
objects, and the paths in $G$ as arrows.  A path is just a map of graphs from a
linear graph into $G$.
Let $\Delta_0$ denote the full subcategory of $\PrSh(\G)$ consisting of the linear graphs.
This category can also be seen as a subcategory of $\Delta$: it has the same
objects but contains only the successor-preserving maps (i.e.~those 
that satisfy $\phi(i+1) = \phi(i)+1$).  In fact the category $\Delta$ can
conveniently be described in terms of $\Delta_0$: it appears by factoring
the composite functor $\Delta_0 \to \PrSh(\G) \to \Cat$ as an 
identity-on-objects functor $j$ followed by a fully faithful functor $i$:
\begin{equation}\label{diagramforDelta}
  \xymatrixrowsep{48pt}
  \xymatrixcolsep{64pt}
  \vcenter{\hbox{\xymatrix @!0
  { &\Delta \ar[r]^i & \Cat \ar@<1ex>[d]^{\text{forgetful}} \ar@{}[d]|{\dashv} \\
\G \ar[r] &\Delta_0 \ar[u]^j \ar[r] & \PrSh(\G) 
\ar@<1ex>[u]^{\text{free}}
}}}
\end{equation}
Among the maps in
$\Delta$ not in $\Delta_0$ are the end-point-preserving 
maps; it is easy to check that every map in $\Delta$ factors uniquely as
an end-point-preserving map followed by a map in $\Delta_0$ (i.e.~a free map).
This factorisation system, a special case of {\em generic/free factorisation} 
\cite{Weber:TAC13}, is an important ingredient in the (modern)
proof of the nerve theorem.

The category $\Delta_0$ has a Grothendieck topology (\cite{MacLane-Moerdijk})
in which a family of maps is
declared to form a cover if they are jointly surjective (on dots as well as on 
strings).  To say that a
simplicial set $X: \Delta\op\to\Set$ satisfies the Segal condition amounts
to saying that its restriction to $\Delta_0$ is a sheaf. 

The generic part of $\Delta$ parametrises the algebraic structure: composition
and identity arrows.  On the other hand, $\Delta_0$ serves to take care of
source-target bookkeeping, and to express the Segal condition.  In the
one-object situation (where $X_0$ is singleton), there is no bookkeeping, and
indeed the notion of monoid can be described solely in terms of the generic part
of $\Delta$.  In fact, the opposite of the category of generic maps in $\Delta$
is the free monoidal category on a monoid, also known as the algebraist's Delta,
or the monoidal Delta, as described in \cite{MacLane:categories}, Ch.VII.

\section{Feynman graphs}

There are various ways to formalise the notion of graph with
open edges (e.g.~\cite{Borisov-Manin:0609748}).
Most of them do not naturally lead to a sensible
notion of morphism.  Although the following definition is very
natural, it seems to be new:

A {\em Feynman graph} is a diagram of finite sets
$$
\xymatrix {
E \ar@(lu,ld)[]_i & H \ar[l]_s \ar[r]^t & V
}
$$
such that $s$ is injective and $i$ is a fixpoint-free involution.  For the
present purposes we also need to impose a connectedness condition.  The set $V$
is the set of {\em vertices}.
The set $H$ is the set of {\em half-edges} or {\em flags}: these are pairs
consisting of a vertex together with the germ of an emanating edge.  Finally the
set $E$ is the set of oriented edges.  The involution $i$ reverses the orientation.
The map $t$ forgets the emanating edge.  The map $s$ returns the emanating edge
in the direction pointing away from the vertex.  A {\em port} is by definition
an (oriented) edge in the complement of the image of $s$.  The set of ports of a 
graph is called its {\em interface}.
An {\em inner edge} is an $i$-orbit both of whose elements are in 
the image of $s$.

From now on we just say {\em graph} for (connected) Feynman graph.
We define the category $\Gr_0$ by taking its objects to be the 
graphs, and its
morphisms to be the diagrams
$$
\xymatrix {
E' \ar@(lu,ld)[] \ar[d] & H' \ar[l] \ar[r] \ar[d]  \drpullback& V' \ar[d] \\
E \ar@(lu,ld)[] & H \ar[l] \ar[r] & V ,
}
$$
the right-hand square being a pullback.
The pullback condition says that each vertex must map to a vertex of
the same valence.  In geometric terms the morphisms are the precisely the
{\em etale} maps (i.e.~local isomorphisms).  (The category $\Gr_0$ will play the 
role $\Delta_0$ plays for categories, as in \S\ref{Sec:classical}.)

An {\em elementary graph} is a graph without inner edges.
Here are the first few elementary graphs:
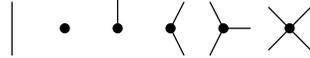
\begin{center}\begin{texdraw}
\lvec (0 20)
\move (20 10) \onedot
\move (40 10) \onedot \lvec (40 21)
\move (60 10) \onedot \lvec (65 20) \move (60 10) \lvec (65 0)
\move (80 10) \onedot \lvec (75 20) \move (80 10) \lvec (75 0) \move (80 10) 
\lvec (90 10)
\move (105 10) \onedot \move (97 18) \lvec (113 2)
\move (97 2) \lvec (113 18)
\end{texdraw}\end{center}
The first one, called the {\em trivial graph} and denoted $\triv$, is given by
$$
\xymatrix {
2 \ar@(lu,ld)[] & 0 \ar[l]\ar[r] & 0 .
}
$$
The remaining ones are of the form 
$$
\xymatrix {
2n \ar@(lu,ld)[] & n \ar[l]\ar[r] & 1 ,
}
$$
for each finite set $n$;  we denote it $n$ again. 
Let $\elGr$ denote the full subcategory of
$\Gr_0$ consisting of the elementary graphs.  We have:
\begin{align*}
  \Hom(\triv,\triv)\ = & \ 2\\
  \Hom(m, n) \ = & \ \begin{cases} n! & \text{ if }m=n\\ 0 & \text{ if } 
  m\neq n\end{cases} \\
  \Hom(\triv , n) \  = & \ 2n \\
  \Hom(n, \triv) \ = & \ 0. 
\end{align*}

It is easy to check that every graph $G$ is canonically a colimit in $\Gr_0$
of its elementary subgraphs.
%
%
%
A family of maps with codomain $G$ is called a {\em cover} of $G$ if it is jointly 
surjective on edges and vertices;  this defines the {\em etale} topology on the
category $\Gr_0$.  The canonical colimit decomposition of a graph is also a 
canonical cover, and it follows readily that there is an 
equivalence of categories between presheaves on $\elGr$ and sheaves on $\Gr_0$:
$$
\PrSh(\elGr) \isopil \Sh(\Gr_0) .
$$
(More formally, the full inclusion of categories
$\elGr \into \Gr_0$ induces an essential geometric embedding of presheaf toposes
$\PrSh(\elGr) \to \PrSh(\Gr_0)$ and it is well-known (see 
\cite{MacLane-Moerdijk}, Ch.VII) that every such induces a unique topology
on $\Gr_0$ giving the above equivalence.)

\section{Graphical species}

A presheaf $F:\elGr\op\to\Set$ is called a {\em graphical species}; its value on $n$ 
is denoted $F[n]$.  Explicitly, a graphical species is
given by an involutive set $C = F[\triv]$, and for each $n\in \N$ a
set $F[n]$ with $2n$ projections to $C$, permuted by a $\mathfrak S_n$-action
on $F[n]$ and by the involution on $C$.  If $C$ is singleton, the classical
notion of species \cite{Joyal:foncteurs-analytiques}, \cite{Bergeron-Labelle-Leroux}
results.

Graphical species parametrise the possible ways of imposing local structure and
decoration on graphs. For each graphical species $F$, the category of
$F$-structured graphs is the comma category $\Gr_0 \comma F$ (i.e.~the category 
whose objects are graphs $G$ equipped with a morphism $NG \to F$ in 
$\PrSh(\elGr)$, where $NG$ denotes the presheaf $n \mapsto \Hom_{\Gr_0}(n,G)$;
the arrows in $\Gr_0 \comma F$ are etale graph maps $G \to G'$ compatible with
the morphisms to $F$).

As examples, there are graphical species for directed graphs, bipartite graphs,
ribbon graphs, and so on.  (In contrast, non-local notions like `graphs of genus
$g$' are not given by graphical species.)
%
Every quantum field theory (see for example
Itzykson-Zuber~\cite{Itzykson-Zuber:QFT}) provides an example of a graphical
species: $F[\triv]$ is the set of field labels, and $F[n]$ is the set of
interaction labels of valence $n$ (multiplied by $\mathfrak
S_n/\Aut$, where $\Aut$ is the symmetry group of the interaction).  For example,
in quantum electrodynamics, there are three field labels, $F[\triv]
=\{e\-,e^+,\gamma\}$ (with the involution intechanging $e^-$ and $e^+$ while
leaving $\gamma$ fixed), and one interaction label
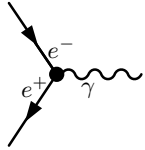
\begin{center}\begin{texdraw}
  \linewd 0.9
  \arrowheadtype t:F \arrowheadsize l:8 w:5
  \move (-20 27) \avec (-9 10.5)
  \move (-10 12) \lvec (-2 0)
  \avec (-14 -18)
  \move (-13 -16.5) \lvec (-20 -27)
  \move (-2 0) \fcir f:0 r:3 \move (0 0) \clvec (1 2.5)(4 2.5)(5 0) \clvec (6 -2.5)(9 -2.5)(10 0) 
  \clvec (11 2.5)(14 2.5)(15 0) \clvec (16 -2.5)(19 -2.5)(20 0)
  \clvec (21 2.5)(24 2.5)(25 0) \clvec (26 -2.5)(29 -2.5)(30 0) 
  \htext (0 10) {\footnotesize $e^-$}
  \htext (-10 -4.5) {\footnotesize $e^+$}
  \htext (10 -6) {\footnotesize $\gamma$}
  \linewd 0.5 
\end{texdraw}\end{center}
Each such label can be applied in $3!$ different ways to a given trivalent 
vertex,  hence $F[3] = 3!$ (and $F[n]=0$ for $n\neq 3$).

\bigskip

The equivalence $\PrSh(\elGr) \simeq \Sh(\Gr_0)$ means that every graphical 
species can be evaluated not only on elementary graphs but on all graphs: if 
$F$ is a graphical species and $G$ is a graph, then
$$
F[G] = \lim_{E\in\elGr \downarrow G} F[E]
$$
where $E$ runs over the category of elements of $G$, i.e.~all the elementary
subgraphs of $G$ and the way they are glued together to give $G$.

\section{Compact symmetry multicategories}

We shall define compact symmetric multicategories as algebras for a monad
defined in terms of sums over graphs.  For the notions of monad and algebras for
a monad we refer to Mac Lane~\cite{MacLane:categories}, Ch.VI.

Let $n$ be a finite set.  
An {\em $n$-graph} is a graph whose set of ports is $n$.
A morphism of $n$-graphs is an isomorphism leaving the set of ports fixed.
We denote this groupoid by $n\kat{-Gr}_\text{iso}$.

We now define the {\em monad for compact symmetric multicategories}:
\begin{align*}
   \PrSh(\elGr) \ \longrightarrow & \ \PrSh(\elGr) \\  
   F \ \longmapsto & \ \ov F ,
\end{align*}
where $\ov F$ is the graphical species given by $\ov F[\triv] :=  
F[\triv]$ and
\begin{align*}
\ov F[n] \ := & \ \colim_{G \in n\kat{-Gr}_\text{iso}} F[G] \\[6pt]
=& \ \sum_{G\in \pi_0(n\kat{-Gr}_\text{iso})}  \frac{F[G]}{\Aut_n(G)} \\[6pt]
=& \ \pi_0 \big( n\kat{-Gr}_\text{iso} \comma F \big) .
\end{align*}
Here the first equation follows since $n\kat{-Gr}_\text{iso}$ is just a 
groupoid: the sum is over isomorphism classes of $n$-graphs, and $\Aut_n(G)$ denotes
the automorphism group of $G$ in $n\kat{-Gr}_\text{iso}$.
The second equation is a lengthy computation with automorphism groups.

This definition is essentially just the coloured version of the definition of
Getzler and Kapranov~\cite{Getzler-Kapranov:9408}.  A formal argument why this
endofunctor has a monad structure can be found in \cite{Getzler-Kapranov:9408}.
Exploiting the third characterisation we can give a heuristic argument (which
can be made into
a formal proof): $\ov F[n]$ is the set of isomorphism classed of
$n$-$F$-graphs: it is the set of ways to decorate $n$-graphs by the graphical
species $F$.  Now $\ov { \ov F} [n]$ is the set of $n$-graphs
decorated by $F$-graphs: this means that each vertex is decorated by a graph
with matching interface.  We can draw each vertex as a circle
with the decorating $F$-graph inside, and the monad structure then consists in
erasing these circles, turning a graph with vertices decorated by $F$-graphs into a
single $F$-graph.

\bigskip

Let $\kat{CSM}$ denote the category of algebras for the monad $F\mapsto \ov F$.
We call its objects {\em compact symmetric multicategories}.  Hence a compact
symmetric multicategory is a graphical species $F:\elGr\op\to\Set$ equipped with
a structure map $\ov F \to F$: it amounts to a rule which for any
$n$-graph $G$ gives a map $F[G] \to F[n]$, i.e.~a way of constructing a single
operation from a whole graph of them.  This rule is subject to a few easy axioms
(cf.~\cite{MacLane:categories}, Ch.VI), amounting roughly to independence of the 
different ways of breaking the computation into steps.

\section{The nerve theorem for compact symmetric multi\-categories}

We now consider the diagram
\begin{equation}\label{diagramforPhi}
  \xymatrixrowsep{48pt}
  \xymatrixcolsep{72pt}
  \vcenter{\hbox{\xymatrix @!0
  {&\Gr \ar[r]^i & \kat{CSM} \ar@<1ex>[d]^{\text{forgetful}} \ar@{}[d]|{\dashv} \\
\elGr \ar[r] &\Gr_0 \ar[u]^j \ar[r] & \PrSh(\elGr) 
\ar@<1ex>[u]^{\text{free}}
}}}
\end{equation}
obtained by factoring $\Gr_0 \to \kat{CSM}$ as an identity-on-objects functor
$j$ followed by a fully faithful functor $i$, just like in
(\ref{diagramforDelta}).  In other words, $\Gr$ is the Kleisli category 
\cite[Ch.VI.5]{MacLane:categories} of the
monad, restricted to $\Gr_0$.  This means that a morphism in $\Gr$ from $G$ to
$G'$ is defined as a morphism of graphical species from $G$ to $\ov{G}{}'$.  So
where the free maps (those coming from $\Gr_0$) send vertices to vertices (of
the same valence), the general maps in $\Gr$ send vertices to `subgraphs' ---
more precisely, a vertex of $G$ is sent to an etale map to $G'$, and the domain
of this etale map must have the same interface as the original vertex.  With
this description is is easy to establish the following factorisation property in
$\Gr$: every map in $\Gr$ factors as a refinement followed by an etale map, in
analogy with the factorisation system in $\Delta$.  The {\em refinements} are
given by taking the domain graph and refining each node, i.e.~replacing the node
by a graph with the same interface.  The etale maps are the free maps, with
respect to the adjunction, while the refinements are the so-called generic maps (in
the weak sense of \cite{Weber:TAC13}), i.e.~characterised by a certain 
universal property.

The embedding $i:\Gr\to\kat{CSM}$ induces the
nerve functor 
\begin{align*}
  N: \kat{CSM} \ &  \longrightarrow \ \PrSh(\Gr)\\
  X \ &\longmapsto \ \Hom_{\kat{CSM}}(i(  \_ ) , X)
\end{align*}
featured in our main theorem:

\bigskip

\noindent \textbf{Theorem.} \cite{Joyal-Kock:FeynmanInPrep} {\em
  The nerve functor $N : \kat{CSM} \to \PrSh(\Gr)$ is fully faithful,
  and a  presheaf is in the essential image of $N$ if and only if it satisfies 
  the Segal condition, i.e.~its restriction to $\Gr_0$ is a sheaf.}

\bigskip

\noindent
The proof follows the ideas and techniques of \cite{Berger:Adv}, 
\cite{Weber:TAC18} and \cite{Kock:0807}.  The main point is to prove that a 
certain left Kan extension is preserved by the monad, which in turn relies on
the generic/free factorisation.  The details will soon be made
available~\cite{Joyal-Kock:FeynmanInPrep}.

\bigskip

One can note, exactly as in the category case, that the generic part encodes the
algebraic structure, whereas the free part is essential for bookkeeping and
expressing the Segal condition.  In the one-object case, 
the category of generic maps is sufficient: its opposite category is essentially
the category of graphs introduced by Getzler and
Kapranov~\cite{Getzler-Kapranov:9408} to study modular operads,
and widely used in the subsequent
literature on the subject.

\nocite{Borisov-Manin:0609748}
\nocite{Moerdijk-Weiss:0701293}

\footnotesize

\hyphenation{mathe-matisk}


\end{document}